\documentclass[12pt]{amsart}

\usepackage{amsmath, amssymb, amsthm}
\usepackage[backend=biber]{biblatex}

\newtheorem{lemma}{Lemma}
\newtheorem{remark}{Remark}
\newtheorem{example}{Example}
\newtheorem{conjecture}{Conjecture}

\title{On the number of even latin squares of even order}
\author{Carolin Hannusch}

\address{Faculty of Informatics, University of Debrecen, Hungary}
\email{hannusch.carolin@inf.unideb.hu}

\begin{document}
	
\maketitle

\begin{abstract}
	We recall the Alon-Tarsi conjecture on the number of even latin squares. We introduce a map which switches the parity of a latin square under certain requirements. An example is included.
\end{abstract}

\section{Preliminaries}

A \emph{latin square} of order $n$ is an $n\times n$ array whose entries are $n$ distinct symbols, each symbol occuring exactly once in each row and each column.
Then each row and column is corresponding to a permutation.

If there is an ordering on the $n$ distinct symbols, then we can determine the parity of each permutation. An \emph{inversion} is a pair of symbols in a permutation such that the order of the pair is not the same as in the natural ordering on the set of $n$ symbols. If the number of inversions in a permutation is even (odd), then the permutation is called even (odd) permutation. We know if a permutation has odd cycle length, then it is an even permutation and if its cycle length is even, then the permutation is odd. The product of two odd permutations or two even permutations is an even permutation, the product of an odd and an even permutation in any order is an odd permutation \cite{jacobson2012basic}. In the current paper, we use the numbers from $1$ to $n$ as distinct symbols with their natural ordering $1<2<\ldots n.$

The Alon-Tarsi conjecture states that for even $n$ the number of even latin squares is not equal to the number of odd latin squares \cite{friedman2019alon}.

For $n\in \{2,4,6,8\},$ we know that the number of even latin squares (A114628 in \cite{sloane2018line}) is greater than the number of odd latin squares of order $n$ (A114629 in \cite{sloane2018line}). 

There exist bounds for the number of latin squares of order $n$ (\cite{colbourn2006latin}, Theorem 1.19), for the number of symmetric latin squares of order $n$ (\cite{ye2011number}) and for the number of reduced latin squares.  In \cite{drisko1997number} it was shown that $Els(n) \neq Ols(n)$ for order $p+1,$ where $p$ is prime. 

\section{Notations}

Let $\lambda$ be a latin square. We denote the row permutations of $\lambda$ by $\sigma_1,\ldots, \sigma_n$ and the column permutations of $\lambda$ by $\tau_1,\ldots, \tau_n.$ A reduced latin square is a latin square with $\sigma_1 = 1_{id}$ and $\tau_1 = 1_{id},$ where $1_{id}$ denotes the identity permutation.

\medskip

Throughout the paper, we will use the following notations:

\begin{center}
\begin{tabular}{|c|c|}
	\hline
	$R_{even}$ & the set of reduced even latin squares\\
	\hline
	$R_{odd}$ & the set of reduced odd latin squares\\
	\hline
	$L_{even}$ & the set of even latin squares\\
	\hline
	$L_{odd}$ & the set of odd latin squares\\
	\hline
	$Erls(n)$ & the number of reduced even latin squares of order $n$\\
	\hline
	$Orls(n)$ & the number of reduced odd latin squares of order $n$\\
	\hline
	$Els(n)$ & the number of even latin squares of order $n$\\
	\hline
	$Ols(n)$ & the number of odd latin squares of order $n$\\
	\hline
\end{tabular}
\end{center}

It is a well-known fact that $Els(n) = Ols(n)$ for odd integers $n,$ since switching two columns in a latin square of odd order changes its sign. If $n$ is even, then this question is more complicated to answer and has been unsolved for a quite long time in general \cite{alon1992colorings} and has been proven for some special cases \cite{drisko1997number}, \cite{drisko1998proof}.

\section{A parity-switching map between latin squares}

In this section, we introduce a map on latin squares, which changes parity under certain circumstances. We show that for each reduced odd latin square there exists at least one map such that its image is an even latin square of the same order. Therefore this map could help to prove that $Els(n)>Ols(n)$ for odd $n.$ 
 Let $\lambda$ be a reduced latin square and $i,j$ the indices of two distinct rows, $i,j\neq 1$. We define $m_{i,j}(\lambda)$ in the following way:
 
 \begin{enumerate}
 	\item Cut off the first row and first column of $\lambda.$ We get an $(n-1)\times (n-1)$ array $\lambda^{\star}$
 	\item Interchange all elements of the rows $i$ and $j$ in $\lambda^{\star}$, except $i$ in row $j$ and $j$ in row $i$ and those symbols, which are in the same cycle as $i$ and $j$ in $\sigma_i$ and $\sigma_j.$
 	\item Add the first row and first column of $\lambda$ to $m_{i,j}(\lambda^{\star})$ and the resulting $n\times n$ array is $m_{i,j}(\lambda).$
 \end{enumerate}

Thus if the rows $i$ and $j$ in $\lambda$ can be written as
$$\begin{array}{cccccc}
	i & : & l_1 & l_2 & \ldots & l_{n-1}\\
	j & : & m_1 & m_2 & \ldots & m_{n-1}
\end{array}$$

and we assume that $l_s, l_t$ are in one cycle with $i$ and $j$ in both row permutations, then $i$ and $j$ in $m_{i,j}(\lambda)$ can be written as

$$\begin{array}{cccccccccccc}

	i & : & m_1 & \ldots & l_s & \ldots & j & \ldots & \ldots & l_t & \ldots & m_{n-1}\\
		j & : & l_1 & \ldots & m_s & \ldots & \ldots & i & \ldots & m_t & \ldots & l_{n-1}\\
\end{array}$$

\begin{remark}
	The map $m_{i,j}$ can be defined similarly on columns of a latin square.
\end{remark}

\begin{example}
	Let $\lambda$ be the following latin square:
	$$\begin{array}{|c|c|c|c|c|c|}
		\hline
		1 & 2 & 3 & 4 & 5 & 6 \\
		\hline
		2 & 1 & 5 & 6 & 3 & 4\\
		\hline
		3 & 6 & 1 & 5 & 4 & 2\\
		\hline
		4 & 3 & 2 & 1 & 6 & 5\\
		\hline
		5 & 4 & 6 & 2 & 1 & 3\\
		\hline
		6 & 5 & 4 & 3 & 2 & 1\\
		\hline
	\end{array}$$

The sign of the product of all row permutations and all column permutations is negative, therefore $\lambda$ is an odd reduced latin square. We have $\sigma_3\sigma_4 = (1,2,5)(3,4,6).$ Thus $m_{3,4}(\lambda)$ is the following latin square:

	$$\begin{array}{|c|c|c|c|c|c|}
	\hline
	1 & 2 & 3 & 4 & 5 & 6 \\
	\hline
	2 & 1 & 5 & 6 & 3 & 4\\
	\hline
	3 & 6 & 2 & 1 & 4 & 5\\
	\hline
	4 & 3 & 1 & 5 & 6 & 2\\
	\hline
	5 & 4 & 6 & 2 & 1 & 3\\
	\hline
	6 & 5 & 4 & 3 & 2 & 1\\
	\hline
\end{array}$$

The sign of the product of all row permutations and all column permutations in this latin square is positive, therefore $m_{3,4}(\lambda)$ is an even reduced latin square.

\end{example}

\begin{lemma}\label{reduced}
	The map $m_{i,j}$ maps a reduced latin square to a reduced latin square for every pair $(i,j)\in \{1,\ldots,n\}\times \{1,\ldots,n\}.$
\end{lemma}

\begin{proof} 

Let $\lambda$ be a reduced latin square. Then $m_{i,j}(\lambda)$ is a latin square, since the rows $i$ and $j$ contain all different $n$ symbols, the other rows do not change. The first column of $m_{i,j}(\lambda)$ is the same as in $\lambda.$ The other columns contain all $n$ different symbols, as they contain the same symbols as the columns of $\lambda,$ only in other order. Further, $m_{i,j}(\lambda)$ is reduced, since its first row and column are the same as in $\lambda$ and $\lambda$ was reduced by assumption. 
\end{proof}

 In Step 2 of the map, it can happen that no element will change. In this case the map $m_{i,j}$ is the identity map. 
 The following lemma helps us to find a condition for $i$ and $j,$ such that $m_{i,j}$ is not the identity map.
 
  \begin{lemma}\label{notidentity}
 	Let $\sigma_i$ and $\sigma_j$ be two row permutations such that $\sigma_i\sigma_j$ is not an $n$-cycle. If $i$ and $j$ are in one cylce of $\sigma_i$ and in one cycle of $\sigma_j,$ then $m_{i,j}$ is not the identity map.
 \end{lemma}
 
 \begin{proof}
 	If $i$ and $j$ are in one cycle in both row permutations, then the map $m_{i,j}$ will fix the elements $i,j$ and the other elements contained in their cycles. Since the product of $\sigma_i$ and $\sigma_j$ is not an $n$-cycle, there are elements in $\{1,\ldots,n\},$ which will be interchanged, therefore $m_{i,j}$ is not the identitdy map.
 \end{proof}

 \begin{lemma}\label{cycle}
 	Let $\sigma_i$ and $\sigma_j$ be two permutations without fixed points. The existence of a cycle, which contains $i$ and $j$ in both, $\sigma_i$ and $\sigma_j$ is equivalent to the fact, that $\sigma_i \sigma_j$ is not an $n$-cycle. 
 \end{lemma}
 
 \begin{proof} 
 First, we assume that $i$ and $j$ are both contained in one cycle of $\sigma_i$ and in one cycle of $\sigma_j.$ Then $\sigma_i\sigma_j$ is not an $n$-cycle. In the other direction, let us assume that $\sigma_i\sigma_j$ is not an $n$-cycle. Then there exists $H\subset \{1,\ldots,n\}$ such that all elements of $H$ are in one cycle in $\sigma_i$ and in one cycle in $\sigma_j.$ Let $H_1$ be the set of elements in the cycle of $i$ in $\sigma_i$ and $H_2$ be the set of elements in the cycle of $j$ in $\sigma_j.$ Since $i\neq j,$ we have that $H_1 \cup H_2 \subset \{1,\ldots,n\}.$ Therefore, we can define the map $m_{i,j}$ on the elements of $H_1 \cup H_2$ and then $m_{i,j}$ is not the identity map.
\end{proof}
 
 \begin{lemma}\label{existence}
 	Let $\lambda \in R_{odd}$ be a reduced odd latin square of order $n,$ where $n$ is an even integer. Then there exist $i,j \in \{1,\ldots,n\},$ such that $m_{i,j}$ is not the identity. 
 \end{lemma}

\begin{proof} 
	Let $\lambda$ be an odd reduced latin square of even order. Then we may assume without loss of generality, that $\prod \tau_i$ is even and $\prod \sigma_i$ is odd (since the transposed of an odd latin square is also an odd latin square). 
Then among $\sigma_1 = 1_{id},\ldots, \sigma_n$ the number of  of odd permutations is odd. Thus among $\sigma_2,\ldots, \sigma_n$ the number of even permutations is even. If this number is $0,$ then there  are $n-1$ odd permutations. The product of any two of them is an even permutation, thus it cannot be an $n$-cycle.
In all other cases, there exist at least two even permutations among $\sigma_2,\ldots, \sigma_n$ such that the product of these is again even, so it cannot be an $n$-cycle. Thus, by Lemma \ref{cycle} and Lemma \ref{existence} the existence of a non-identity $m_{i,j}$ is proven. 
\end{proof}

\begin{lemma}\label{sign}
	Let $\lambda \in R_{odd}$ and $m_{i,j} \neq 1_{id}.$ Then $m_{i,j}(\lambda) \in R_{even}.$
\end{lemma}

\begin{conjecture}\label{number1}
	Let $n\in \mathbb{N}.$ Then $Erls(n) \geq Orls(n).$
\end{conjecture}

 It is well known that $(Els(n)+Ols(n)) = n!(n-1)!(Erls(n)+Orls(n)).$ Since $n!(n-1)!$ is the number of all possible permutations of rows and columns and each permutation can be represented by a finite number of transpositions (i.e.~switching two rows or two columns, respectively) \cite{clark1984elements}. It is also well known that switching two rows (or two columns) does not change the sign of an even latin square of even order. Therefore, if Conjecture \ref{number1} is true, then $Els(n) \geq Ols(n)$ follows.

\subsection{Even latin squares}
On the other way round, applying the map $m_{i,j}$ to an even reduced latin square of even order can result in an even reduced latin square for all possible maps $m_{i,j}.$ See the following example.
\newline $\lambda:$
$$\begin{array}{|c|c|c|c|c|c|c|c|}
	\hline
	1 & 2 & 3 & 4 & 5 & 6 & 7 & 8\\
	\hline
	2 & 3 & 4 & 5 & 6 & 7 & 8 & 1\\
	\hline
	3 & 4 & 5 & 6 & 7 & 8 & 1 & 2\\
	\hline
	4 & 5 & 6 & 7 & 8 & 1 & 2 & 3\\
	\hline
	5 & 6 & 7 & 8 & 1 & 2 & 3 & 4\\
	\hline
	6 & 7 & 8 & 1 & 2 & 3 & 4 & 5\\
	\hline
	7 & 8 & 1 & 2 & 3 & 4 & 5 & 6\\
	\hline
	8 & 1 & 2 & 3 & 4 & 5 & 6 & 7\\
	\hline
	\end{array}$$

\bigskip

For all pairs $(i,j)\in \{2,\ldots,8\} \times \{2,\ldots,8\}$ with $i\neq j$ one of the following three cases occur:
\begin{itemize}
	\item $m_{i,j}$ is the identity for $(i,j)\in \{(2,3), (2,5), (2,7), (3,4), (3,6), \newline (3,8), (4,5), (4,7), (5,6), (5,8), (6,7), (7,8)\}$
	\item $m_{i,j}$ switches the elements of a $4$-cycle and therefore the signs of the product of the column permutations and of the row permutations do not change, this is the case for the pairs $(i,j)\in \{(2,4), (2,8), (3,5), (4,6), (5,7), (6,8)\}$ or 
	\item $m_{i,j}$ switches the elements of a transposition and therefore the signs of the product of the column permutations and of the row permutations do not change again for $(i,j)\in \{(2,6), (3,7), (4,8)\}.$
\end{itemize} 

\printbibliography

\end{document}